\documentclass[a4paper,11pt]{article}

\usepackage{latexsym}
\usepackage{amsmath,amsfonts,amssymb,amsthm}
\newtheorem{thm}{Theorem}[section]

\newtheorem{cor}[thm]{Corollary}
\newtheorem{prop}[thm]{Proposition}

\newtheorem{rem}[thm]{Remark}

\input epsf

\title{On a family of tridiagonal matrices}
\author{Roland Bacher}
\begin{document}
\maketitle

{\sl Abstract: 
We show that certain integral positive definite 
symmetric tridiagonal matrices of determinant $n$ are in 
one to one correspondence with
elements of $(\mathbb Z/n\mathbb Z)^*$. 
We study some properties of this correspondence.

In a somewhat unrelated second part we discuss a construction
which associates a sequence of integral polytopes to
every integral symmetric matrix.}

\section{Introduction}

Let $\alpha_1,\dots,\alpha_n$ be a finite sequence of elements in a
commutative ring.
We denote by $T(\alpha_1,\dots,\alpha_n)$ the symmetric tridiagonal matrix
$$\left(\begin{array}{cccccccccccccccccc}
\alpha_1&1&0&\dots&0&0&0\\
1&\alpha_2&1&&0&0&0\\
0&1&\alpha_3&&0&0&0\\
\vdots&&&\ddots&&&\vdots\\
0&0&0&&\alpha_{n-2}&1&0\\
0&0&0&&1&\alpha_{n-1}&1\\
0&0&0&\dots&0&1&\alpha_n\end{array}\right)$$
of size $n\times n$ with diagonal entries $\alpha_1,\dots,\alpha_n$
and sub and super diagonals consisting only of ones.

Leighton and Newman discovered in \cite{LN} the remarkable fact
that the number of such matrices of size $n\times n$ which are
unimodular, integral and positive definite is given by the 
$n-$th Catalan number
${2n\choose n}/(1+n)$, see also Exercice 6.19 nnn in \cite{St2} and
\cite{Sh}. One should mention that this result was already 
implicitely contained in \cite{CC}, see the solution to 
problem (18).

In the present paper, we are interested in the number of such 
matrices of arbitrary size and determinant $N$ which are integral and 
have diagonal entries in $\{2,3,4,\dots\}$. Such matrices are always
positive definite, the case $\alpha_1=\alpha_2=\dots=\alpha_n=2$ 
defining the root lattice of type $A_n$. 

\begin{thm} \label{thmmain} The number of finite symmetric tridiagonal 
integral matrices of the form 
$T(\alpha_1,\alpha_2,\dots)$ with diagonal entries $\geq 2$ and
determinant $N$ is given by Euler's totient function 
$$\Phi(N)=N\prod_{p\vert N}\left(1-\frac{1}{p}\right)$$
(where the product is over all prime divisors of $N$)
counting the number $\sharp((\mathbb Z/N\mathbb Z)^*)$
of invertible elements of the ring $\mathbb Z/N\mathbb Z$.
\end{thm}

The proof is bijective: Given $(a,N)$ for $a$ in
$\{1,2,\dots,N-1\}$ invertible modulo $N$, we construct 
a finite sequence $(\alpha_1,\alpha_2,\dots)$ giving rise 
to a matrix $T(\alpha_1,\alpha_2,\dots)$ with the desired properties.
This construction is based on properties of $\mathrm{SL}_2(\mathbb Z)$ 
related to Euclid's algorithm.

The organisation of the paper is as follows:

The next section exposes a few useful (and probably well-known) facts
concerning $\mathrm{SL}_2(\mathbb Z)$.

The proof of Theorem \ref{thmmain} is contained in 
Section \ref{sectproof}.

Section \ref{sectalgor} discusses some algorithmic aspects 
of the construction wich associates to an element $a\in
(\mathbb Z/N\mathbb Z)^*$ the corresponding tridiagonal 
matrix of the form $T(\alpha_1,\dots)$. 

Section \ref{sectsigma} describes properties
of the function $a\longmapsto \mathop{tr}(T)-\dim(T)$
for $a\in(\mathbb Z/N\mathbb Z)^*$.

All matrices described by Theorem \ref{thmmain}
are definite positive. Section \ref{sectlattices} 
discusses briefly Euclidean lattices 
admitting the matrices of Theorem \ref{thmmain} as Gram matrices.

The final Section \ref{sectpolyt} is not closely related to
the rest of the paper. It describes a construction of a sequence
of integral polytopes associated to a symmetric integral matrix
and discusses some easy properties of this construction.
The first polytope obtained by applying this construction 
to tridiagonal matrices with $1'$s on the sub and super 
diagonal end zeros elsewhere yields for instance integral 
polytopes with vertices enumerated by Catalan numbers. These polytopes
are however distinct from associahedrons (or Stasheff polytopes).

\section{$\mathrm{SL}_2$ and the determinant of $T(\alpha_1,\dots,
\alpha_n)$}

We consider the map
$$\alpha\longmapsto M(\alpha)=\left(\begin{array}{cc}0&-1\\1&\alpha
\end{array}\right)$$
extended to finite sequences by setting
$$M(\alpha_1,\alpha_2,\dots,\alpha_n)=M(\alpha_1)M(\alpha_2)
\cdots M(\alpha_n)\ .$$
Denoting by
$$\vert T(\alpha_1,\dots,\alpha_n)\vert=\det(T(\alpha_1,\dots,\alpha_n))$$
the determinant of $T(\alpha_1,\dots,\alpha_n)$ we have the following result:

\begin{prop} \label{propdetident}
We have for all $n\geq 2$ the identity
$$M(\alpha_1,\dots,\alpha_n)
=\left(\begin{array}{cc}
-\vert T(\alpha_2,\dots,\alpha_{n-1})\vert&-\vert T(\alpha_2,\dots,
\alpha_n)\vert\\
\vert T(\alpha_1,\dots,\alpha_{n-1})\vert&\vert T(\alpha_1,\dots,
\alpha_n)\vert
\end{array}\right)$$
where we use the convention 
$\vert T(\alpha_2,\dots, \alpha_{n-1})\vert=1$ for $n=2$.
\end{prop}

The proof is a straightforward computation left to the reader.

Proposition \ref{propdetident} implies easily the identity
$$\vert T((\alpha_1,\dots,\alpha_{n+1})\vert\ 
\vert T(\alpha_2,\dots,\alpha_{n-1})\vert
=\vert T(\alpha_1,\dots,\alpha_n)\vert\ 
\vert T(\alpha_2,\dots,\alpha_{n+1})\vert-1$$
which is a particular case of Dodgon's condensation formula.

Another useful result is the following ``symmetry''.

\begin{prop} \label{propab}
If 
$$M(\alpha_1,\alpha_2,\dots,\alpha_n)=\left(\begin{array}{cc}
-a&-b\\c&d\end{array}\right)$$
then 
$$M(\alpha_n,\dots,\alpha_2,\alpha_1)=
\left(\begin{array}{cc}
-a&-c\\b&d\end{array}\right)\ .$$
\end{prop}

{\bf Proof} Introducing the involution 
$V=\left(\begin{array}{cc}0&1\\1&0\end{array}\right)$,
a short computation shows $VM(\alpha)V=M(\alpha)^{-1}$.
This yields
$$VM(\alpha_1,\dots,\alpha_n)^{-1}V=M(\alpha_n,\dots,\alpha_1)$$
which establishes the result by computing
$$\left(\begin{array}{cc}0&1\\1&0\end{array}\right)
\left(\begin{array}{cc}d&b\\-c&-a\end{array}\right)
\left(\begin{array}{cc}0&1\\1&0\end{array}\right)\ .$$
\hfill$\Box$

The easy identity 
$$M(x)M(y)=M(x+1)M(1)M(y+1)=\left(\begin{array}{cc}
-1&-y\\x&xy-1\end{array}\right)$$
shows that the set of integral positive definite matrices 
of the form $T(\alpha_1,\alpha_2,\dots)$ with given determinant is 
infinite without further conditions. 

Leighton and Newman consider matrices of
the form $T(\alpha_1,\dots,\alpha_n)$ which are integral, positive definite 
and have fixed size $n\times n$ in order to get finite sets.
In the present paper, the size of the matrices $T(\alpha_1,\alpha_2,\dots)$
is arbitrary but the diagonal
coefficients $\alpha_1,\alpha_2,\dots$ are restricted to the set 
$\{2,3,4,\dots\}$ of natural integers $\geq 2$. This restriction is
motivated by the following result. 

\begin{prop} \label{propineqfond}
Let $\alpha_1,\alpha_2,\dots,\alpha_l$ be a sequence
of $l\geq 1$ real numbers $a_i\geq 2$ for $i=1,\dots,l$. Consider
$$M(\alpha_1,\alpha_2,\cdots,\alpha_l)=
\left(\begin{array}{cc}-a&-b\\c&d\end{array}\right)\ .$$
We have then the inequalities
$$0\leq a<\min(b,c),\ \max(b,c)<d\hbox{ and }c-a\leq d-b\ .$$
\end{prop}

\begin{cor} \label{cordetcroiss} We have
$$\det(T(\alpha_1,\alpha_2,\dots,\alpha_{l-1}))<\det
(T(\alpha_1,\alpha_2,\dots,\alpha_{l-1},\alpha_l))$$
if $\alpha_1,\dots,\alpha_l$ are real numbers $\geq 2$.

In particular, the number of matrices of the form 
$T(\alpha_1,\alpha_2,\dots)$
with determinant $N$ and integral diagonal coefficients 
$\alpha_1,\alpha_2,\dots\subset\{2,3,4,\dots\}$
is finite.
\end{cor}

{\bf Proof} The first part follows by applying Proposition 
\ref{propdetident} and the inequality $c<d$ 
of Proposition \ref{propineqfond} to the matrix 
$M(\alpha_1,\dots,\alpha_l)$.

The second part follows also from Proposition 
\ref{propdetident} and from the observation that the number of 
integral matrices satisfying the inequalities of
Proposition \ref{propineqfond} with $d=N$ is finite.
\hfill$\Box$

\begin{cor} \label{cordefpos} A matrix of the form 
$T(\alpha_1,\alpha_2,\dots,\alpha_l)$ with
$\alpha_i\geq 2$ for $i=1,2,\dots,l$ is positive definite.
\end{cor}

\begin{rem} Corollary 
\ref{cordefpos} follows of course also from the well-known observation 
that $T(2,2,\dots,2)$ is a Gram matrix for a root lattice of type 
$A$ and from the fact that real symmetric positive matrices 
form a convex cone.
\end{rem}

{\bf Proof of Proposition \ref{propineqfond}} 
By induction on $l$. The results hold obviously
for $l=1$ since
$$M(\alpha_1)=\left(\begin{array}{cc}0&-1\\1&\alpha_1\end{array}\right)$$
satisfies all inequalites for $\alpha_1\geq 2$.

For the induction step we consider
$$\left(\begin{array}{cc}-a&-b\\c&d\end{array}\right)
\left(\begin{array}{cc}0&-1\\1&x\end{array}\right)=
\left(\begin{array}{cc}-b&a-xb\\d&-c+xd\end{array}\right)\ .$$

We have obviously $0<b,d$ and $d>b$ by induction.
Since $0\leq a<b$ we have 
$$-a+xb>(x-1)b\geq b>0$$
for $x\geq 2$.

The inequalities $0<c<d$ show similarly
$$-c+xd>(x-1)d\geq d>0\ .$$ 

The inequalities $c-a\leq d-b$ and $d>b$ imply
$$(xd-c)-(xb-a)=x(d-b)-(c-a)\geq (x-1)(d-b)\geq d-b>0$$
for $x\geq 2$.
\hfill$\Box$

\section{Proof of Theorem \ref{thmmain}}\label{sectproof}

\begin{thm} \label{thmalgo}
For every integral unimodular matrix
$$A=\left(\begin{array}{cc}-a&-b\\c&d\end{array}\right)$$
such that 
$$0\leq \min(a,b,c,d)\hbox{ and }\max(b,c)<d$$
there exists a unique integer $l\geq 1$ and a unique finite sequence 
$(\alpha_1,\dots,\alpha_l)\in\{2,3,4,\dots\}^l$
such that 
$$A=M(\alpha_1)M(\alpha_2)\cdots M(\alpha_l)\ .$$
\end{thm}

Theorem \ref{thmmain} is now implied by the following result.

\begin{prop} \label{propconstr}
Integral matrices of the form $T(\alpha_1,\dots,\alpha_n)$
with determinant $N$ and diagonal coefficients $\geq 2$
are in bijection with the subset of integers in
$\{1,\dots,N\}$ which are invertible modulo $N$.
\end{prop}

{\bf Proof of Proposition \ref{propconstr}} 
Given a pair $(a,N)$ of natural integers with $a\in
\{1,\dots,N-1\}$ such that $a$ is invertible modulo $N$, we consider $b\in
\{1,\dots,N-1\}$ such that $ab\equiv 1\pmod N$ and we set $k=
\frac{ab-1}{N}\in\mathbb N$.
The matrix $A=\left(\begin{array}{cc}-k&-b\\a&N\end{array}\right)$
satisfies the hyptheses of Theorem \ref{thmalgo} and has thus a 
unique factorisation of the form 
$$A=M(\alpha_1)M(\alpha_2)\cdots M(\alpha_l)$$for some integer $l\geq 1$ 
and $\alpha_1,\alpha_2,\dots\in\{2,3,\dots\}^l$.

Proposition \ref{propdetident} shows thus that the matrix 
$T(\alpha_1,\dots,\alpha_l)$ has determinant $N$.

This construction yields an injective map from 
$(\mathbb Z/N\mathbb Z)^*$ into integral matrices of the form
$T(\alpha_1,\alpha_2,\dots)$ with determinant $N$ and diagonal
coefficients $\alpha_i\geq 2$.

The fact that this map is onto follows by applying the unicity result
of Theorem \ref{thmalgo} to the matrix 
$$\left(\begin{array}{cc}
-\vert T(\alpha_2,\dots,\alpha_{n-1})\vert&-\vert T(\alpha_2,\dots,
\alpha_n)\vert\\
\vert T(\alpha_1,\dots,\alpha_{n-1})\vert&\vert T(\alpha_1,\dots,
\alpha_n)\vert
\end{array}\right)=M(\alpha_1)\cdots M(\alpha_l)$$
associated to a suitable tridiagonal matrix 
$T(\alpha_1,\alpha_2,\dots)$.
\hfill$\Box$

{\bf Proof of Theorem \ref{thmalgo}} 
If $a=0$ then unimodularity of $A$ implies $b=c=1$
and $A=M(d)$ for some integer $d\geq 2>1=b=c$. Otherwise, the identity
$$bc=ad+1$$
implies $bc\equiv 1\pmod d$.
This shows that the integers $c,d$ are coprime and
determine $a,b$ uniquely by considering the unique integer
$b\in\{1,\dots,d-1\}$ such that $bc\equiv 1\pmod d$
and by setting $a=\frac{bc-1}{d}\geq 1$.
Unimodularity of $A$ and strict positivity of $a$ imply that
Defining $t\in\{1,\dots,b-1\}$ as the unique integer 
satisfying the congruence $t\equiv -d\pmod b$,
the equality $\det(A)=-ad+bd=1$ implies $at\equiv 1\pmod b$
and we can thus consider the natural integer $s=\frac{at-1}{b}\geq 0$.
The computation
$$\left(\begin{array}{cc}-a&-b\\c&d\end{array}\right)
\left(\begin{array}{cc}b&t\\-a&-s\end{array}\right)=
\left(\begin{array}{cc}0&-1\\ 1&ct-sd\end{array}\right)$$
and the inequality
$$b\left(ct-\frac{at-1}{b}d\right)=bct-adt+d=
t+d>b$$
show thus 
$$\left(\begin{array}{cc}-a&-b\\c&d\end{array}\right)=
\left(\begin{array}{cc}0&-1\\1&ct-sd\end{array}\right)
\left(\begin{array}{cc}-s&-t\\a&b\end{array}\right)$$
where $ct-sd\geq 2$. Since $a=\frac{bc-1}{d}<b$, the integral 
unimodular matrix
$$\left(\begin{array}{cc}-s&-t\\a&b\end{array}\right)$$
satisfies again the hypotheses of Theorem \ref{thmalgo}.
Induction on $d>b$ establishes the existence of a sequence $(\alpha_1,
\alpha_2,\dots,\alpha_l)$ in
$\{2,3,4,\dots\}^l$ such that 
$A=M(\alpha_1)M(\alpha_2)\cdots M(\alpha_l)$.

Unicity of the sequence $\alpha_1,\dots$ is obvious for $a=0$
by Proposition \ref{propdetident}.
For $a>0$ it follows from the observation that the inequalities
$0\leq \min(\tilde s,\tilde t)$ and $\tilde t<b$ determine
the integral unimodular matrix
$$\left(\begin{array}{cc}
-\tilde s&-\tilde t\\a&b\end{array}\right)=\left(\begin{array}{cc}
\tilde\alpha_1&1\\-1&0\end{array}\right)\left(\begin{array}{cc}
-a&-b\\c&d\end{array}\right)=T(\tilde\alpha_2,\tilde\alpha_3,\dots)$$
uniquely. This shows $\tilde\alpha_1=\alpha_1$ and by 
induction $\tilde\alpha_i=\alpha_i$ and $l'=l$ for 
$\tilde\alpha_1,\dots,\tilde\alpha_{l'}\in\{2,3,\dots\}^{l'}$
such that $M(\alpha_1,\dots,\alpha_l)=M(\tilde\alpha_1,\dots,
\tilde\alpha_{l'})$.
\hfill$\Box$

\begin{rem} Theorem \ref{thmalgo} is equivalent to the assertion that
the matrices $M(2),M(3),M(4),\dots$ generate a free submonoid 
of $\mathrm{SL}_2(\mathbb Z)$. The subgroup generated by 
$M(2),M(3),\dots$ is however not free and coincides with $\mathrm{SL}_2(
\mathbb Z)$.
Indeed, the identity
$$M(\lambda)M(\lambda+\mu)^{-1}M(\lambda)=M(\lambda-\mu)$$
shows that the subgroup generated by two integral matrices
$M(k)$ and $M(k+1)$ contains the generators
$$M(0)=\left(\begin{array}{cc}0&-1\\1&0\end{array}\right)\hbox{ and }
M(1)=\left(\begin{array}{cc}0&-1\\1&1\end{array}\right)$$
of $\mathrm{SL}_2(\mathbb Z)$.
\end{rem}

\section{Algorithmic aspects} \label{sectalgor}

In the sequel, $x\pmod y$ denotes always the unique natural integer 
in $\{1,\dots,y\}$ representing the equivalence class of $x$
in $\mathbb Z/y\mathbb Z$ for an integer $x$ and 
a natural integer $y\geq 2$. Similarly, if $x$ is invertible 
modulo $y$, then $x^{-1}\pmod y$ denotes the unique integer
in $\{1,\dots,y-1\}$ such that $x(x^{-1}\pmod y)\equiv 1\pmod y$.

Given a natural integer $N\geq 2$ and an integer $a\in\{1,\dots,N-1\}$ 
which is invertible modulo $N$, we denote by 
$$W(a,N)=\alpha_1,\dots,\alpha_l$$
the unique finite sequence of natural integers $\alpha_i\geq 2$ such that
$$M(\alpha_1\dots\alpha_l)=\left(\begin{array}{cc}
-x&-y\\a&N\end{array}\right)\in \mathrm{SL}_2(\mathbb Z)$$
where $y\equiv a^{-1}\pmod N$ and $x=\frac{ay-1}{N}$,
cf Theorem \ref{thmalgo}.

We denote by $l(a,N)=l$ the length of the word $W(a,N)$
and by $W(a,N)_i=\alpha_i$ the $i-$th element (for $i=1,\dots,l$) of the
sequence $W(a,N)$. We have thus
$$W(a,N)=W(a,N)_1,W(a,N)_2,\dots,W(a,N)_{l(a,N)}\ .$$

We have obviously
$$W(1,N)=W(1,N)_1=N,\ l(1,N)=1\ .$$

If $a,b\in\{1,\dots,N-1\}$ are invertible elements modulo $N$ such that
$ab\equiv 1\pmod N$, Proposition \ref{propab} shows that we have
$$l(a,N)=l(b,N)$$ 
and
$$W(a,N)_i=W(b,N)_{l(a,N)+1-i}\ .$$
Otherwise stated,
the matrices $T(W(a,N))$ and $T(W(b,N))$ are related to each other 
by conjugation with 
the antidiagonal involution 
of size $l(a,N)\times l(a,N)$.

We denote by $\rfloor x\lfloor$ the integral part of a real number $x$
defined as the unique integer satisfying
$$x-1<\lfloor x\rfloor \leq x\ .$$

The following result allows computations:

\begin{thm} \label{thmcoeff} For every integer $N\geq 2$ we have 
$$W(1,N)=N\ .$$

For an integer $a\in\{2,\dots,N-1\}$ invertible modulo $N$
(with $N$ a natural integer $\geq 3$) we have
$$W(a,N)=W\left((-N)\pmod a,a\right),
\left(1+\left\lfloor\frac{N}{a}\right\rfloor\right)$$
or equivalently
$$W(a,N)=\left(1+\left\lfloor\frac{N}{b}\right\rfloor\right),
W\left((-N)^{-1}\pmod b,b\right)$$
where $b=(a^{-1}\pmod N)\in \{1,\dots,N-1\}$.
\end{thm}

\begin{rem} The map $(a,N)\longmapsto ((-N\pmod a),a)$ 
involved in the first equality for $W(a,N)$ for $a\geq 2$ coincides
(up to the $-$ sign) with Euclid's celebrated algorithm
for computing the greatest common divisor 
of $a$ and $N$.
\end{rem}

{\bf Proof of Theorem \ref{thmcoeff}} 
The assertion concerning $W(1,N)$ is obvious.

Set $k=\frac{ab-1}{N}\geq 1$ where $b=a^{-1}\pmod N$ is the 
inverse of $a$ modulo $N$ in $\{2,\dots,N-1\}$. Computing
$$\left(\begin{array}{cc}-k&-b\\a&N\end{array}\right)
\left(\begin{array}{cc}1+\lfloor N/a\rfloor&1\\-1&0\end{array}\right)=
\left(\begin{array}{cc}
-k(1+\lfloor N/a)\rfloor+b&-k\\a(1+\lfloor N/a\rfloor)-N&a\end{array}
\right)\ ,$$
the lower left coefficient of the last matrix is given by 
$$a\left(1+\frac{N-(N\pmod a)}{a}\right)-N=a+N-(N\pmod a)-N=(-N)\pmod a$$
and is thus a strictly positive natural integer $<a$. 
Since 
$$k(a(1+\lfloor N/a\rfloor)-N)\equiv -kN\equiv 1\pmod a\ ,$$
positivity of $k,a(1+\lfloor N/a\rfloor)-N$, the obvious inequality
$k<a$ and unimodularity of
all involved matrices imply that  the last matrix 
satisfies the conditions of Theorem \ref{thmcoeff}.
This implies the first equality for $a\geq 2$ by induction on $N$.

The second equality can be proven similary.
It follows also from the first equality and from Proposition 
\ref{propab}.
\hfill$\Box$

For computing $W(a,N)$ with $a\in\{1,\dots,N-1\}$ invertible modulo $N$
and very close to a huge integer $N$ the following result is useful:

\begin{prop} \label{propalgo}
The sequence $W(N-1,N)$ is the constant sequence $2,2,\dots,2$ of 
length $N-1$.

We have for $a\in\{2,\dots,N-2\}$ invertible modulo $N$ the equality
$$W(a,N)=W((\lambda+1)a-\lambda N,\lambda a-(\lambda-1)N),2,2,2,\dots,2$$
where $\lambda=\lfloor a/(N-a)\rfloor$ and where $2,2,2,\dots,2$
is the constant sequence of length $\lambda$.
\end{prop}

{\bf Proof} The assertion for the sequence $W(N-1,N)$ follows
from the case $k=N-1$ of the identity
$$\left(\begin{array}{cc}0&-1\\1&2\end{array}\right)^k=
\left(\begin{array}{cc}1-k&-k\\k&k+1\end{array}\right)$$
which is easily established.

For $a\geq 2$ such that $a<N/2$ we have $a/(N-a)<1$ yielding $\lambda=0$
and the formula is trivial.

The case $a=N/2$ implies $(a,N)=(1,2)$ corresponding to
the case $(N-1,N)$ for $N=2$.

The proof in the case $a>N/2$ is by induction on $\lambda$.
If $\lambda=1$ we have to prove
$$W(a,N)=W(2a-N,a),2$$
where $\frac{1}{2}N<a<\frac{2}{3}N$.  
This boils down to the identity
$W(a,N)=W((-N)\pmod a,a),(1+\lfloor N/a\rfloor)$
of Theorem \ref{thmcoeff}.

For $\frac{2}{3}N\leq a<N-1$ we have as above
$W(a,N)=W(2a-N,a),2$ by Theorem \ref{thmcoeff}.
The trivial identities
$$\frac{2a-N}{a-(2a-N)}=\frac{a}{N-a}-1$$
and $$k(2a-N)-(k-1)a=(k+1)a-kN$$ 
and induction on $\lambda$ end the proof.\hfill$\Box$

Theorem \ref{thmcoeff} and Proposition \ref{propalgo}
give a fast algorithm for computing 
$$\sum_{i=1}^{l(a,N)}(W(a,N)_i)^e$$
for arbitrary $e\in \mathbb N$ (or $e\in\mathbb C$). The case $e=0$ 
corresponds to
the length $l(a,N)$ of the sequence $W(a,N)$
and the case $e=1$ yields the trace of the tridiagonal matrix 
$T$ associated to $W(a,N)$. The details are as follows:

INPUT: $(a,N)$ with $a\in\{1,\dots,N-1\}$ invertible modulo 
$N$ (for $N\geq 2$) and an integer $e$.

Set $r=0$.

LOOP: 

IF $a=1$, OUTPUT: $(r+N^e)$ and STOP ENDIF.

IF $a<(N/2)$, replace $r$ by $r+(1+\lfloor N/a\rfloor)^e$
and $(a,N)$ by $((-N)\pmod a,a)$ where $(-N)\pmod a$ is
the choosen in $\{1,\dots,a-1\}$) ENDIF.

IF $a=N-1$, OUTPUT: $r+(N-1)2^e$ and STOP ENDIF.

IF $a>N/2$, set $\lambda=\lfloor a/(N-a)\rfloor$, 
replace $r$ by $r+\lambda 2^e$ and $(a,N)$
by $((\lambda+1)a-\lambda N,\lambda a-(\lambda-1) N)$ ENDIF.

END OF LOOP.

\section{A function with an additional symmetry}\label{sectsigma}

Set $\sigma(a,N)=\sum_{k=1}^{l(a,N)}(W(a,N)_k-1)=
\hbox{tr}(T(W(a,N)))-\dim(T(W(a,N)))$ for 
$a\in\{1,\dots,N-1\}$ invertible modulo $N$.

Proposition \ref{propab} shows that the function $\sigma$
satisfies
$\sigma(a,N)=\sigma((a^{-1}\pmod N),N)$.

The following result shows that $\sigma$ satifies an additional symmetry.
 
\begin{thm} \label{thmsymsigma} We have 
$$\sigma(a,N)=\sigma(N-a,N)$$
for every $a\in\{1,\dots,N-1\}$ which is invertible modulo $N$.
\end{thm}

{\bf Proof} We have 
$\sigma(1,N)=N-1$ and $\sigma(N-1,N)=(N-1)(2-1)=N-1$.

Consider now $a$ such that $2\leq a<\frac{N}{2}$. Setting 
$\lambda=\lfloor (N-a)/a\rfloor=\frac{N-(N\pmod a)}{a}-1$, Proposition
\ref{propalgo} shows
$$\sigma(N-a,N)=\sigma(N-(\lambda+1)a,N-\lambda a)+\lambda$$
implying
$$\sigma(N-a,N)=\sigma(a,N-\lambda a)+\lambda$$
by induction on $N$. Theorem \ref{thmcoeff} yields thus
$$\sigma(N-a,N)=\sigma((-N)\pmod a,a)+
\frac{N-\lambda a-(N\pmod a)}{a}+\lambda$$
$$=\sigma((-N)\pmod a,a)
+\lfloor N/a\rfloor$$
which establishes the result since we have also  
$$\sigma(a,N)=\sigma((-N)\pmod a,a)+\lfloor N/a\rfloor$$
by Theorem \ref{thmcoeff}.\hfill$\Box$

A few values of $\sigma(a,N)$
are given by the following Table:
$$\begin{array}{lrrrrrrrrrrrrrrrrrrrrrr}
a=&1&2&3&4&5&6&7&8&9&10&11&12&13&14&15&16\\
N=2:&1\\
N=3:&2&2\\
N=4:&3&&3\\
N=5:&4&3&3&4\\
N=6:&5&&&&5\\
N=7:&6&4&4&4&4&6\\
N=8:&7&&4&&4&&7\\
N=9:&8&5&&5&5&&5&8\\
N=10:&9&&5&&&&5&&9\\
N=11:&10&6&5&5&6&6&5&5&6&10\\
N=12:&11&&&&5&&5&&&&11\\
N=13:&12&7&6&6&5&7&7&5&6&6&7&12\\
N=14:&13&&6&&6&&&&6&&6&&13\\
N=15:&14&8&&6&&&8&8&&&6&&8&14\\
N=16:&15&&7&&7&&6&&6&&7&&7&&15\\
N=17:&16&9&7&7&6&7&6&9&9&6&7&6&7&7&9&16
\end{array}$$

\subsection{Continued fraction expansions}

Let $[\gamma_1,\gamma_2,\dots]$
be the continued fraction expansion of $x\in(0,1)$ defined recursively 
by $x=\frac{1}{\gamma_1 +[\gamma_2,\dots]}$
where $\gamma_1=\lfloor 1/x\rfloor$ and $[\gamma_2,\dots]=1/x-\gamma_1$ 
(and $0=[\emptyset]$ by convention).

\begin{thm} \label{thmconfrac} We have
$$1+\sigma(a,N)=\sum_{i\geq 1} \gamma_i$$
where $\frac{a}{N}=[\gamma_1,\gamma_2,\dots]$
is the continued fraction expansion
of the reduced rational fraction $a/N$ in $(0,1)$.
\end{thm}

{\bf Proof} The result holds obviously for $a=1$.  
For $a>1$, we have
$$\gamma_1=\lfloor N/a\rfloor$$ and
$$[\gamma_2,\dots]=N/a-\lfloor N/a\rfloor=\frac{N-(N-(N\pmod a))}{a}=
\frac{N\pmod a}{a}$$
The equality $\gamma_1=\alpha_1$ and induction on $N$
end the proof.\hfill$\Box$

\begin{rem} Theorem \ref{thmsymsigma} follows also from Theorem 
\ref{thmconfrac}
together with the easy continued fraction expansion
$$1-x=[1,\alpha_1-1,\alpha_2,\alpha_3,\dots]$$
for $x=[\alpha_1,\alpha_2,\alpha_3,\dots]$ in $(0,1/2)$.
\end{rem}

\section{Lattices} \label{sectlattices}

Given $d$ integers $\alpha_1,\dots,\alpha_d\in\{2,3,\dots\}^d$, 
we denote by $\Lambda(\alpha_1,\dots,\alpha_d)$ the Euclideean
lattice of rank $d$ with scalar product given by the Gram matrix
$T(\alpha_1,\dots,\alpha_d)$ (which is positive definite by Corollary 
\ref{cordefpos}. We say that such a lattice is of type T.

\begin{prop} \label{propeucllatt} A lattice $\Lambda(\alpha_1,\dots,\alpha_d)$
of type $T$ has minimal norm $
\min(\alpha_1,\dots,\alpha_d)\geq 2$.
\end{prop}

{\bf Proof} We work with a basis with Gram 
matrix $T(\alpha_1,\dots,\alpha_d)$.
Replacing $\alpha_1,\dots,\alpha_d$ by $2$ yields a root lattice
$\Lambda'$ of type $A_d$ and the obvious map 
sending a vector $\lambda=(\lambda_1,\dots,\lambda_d)$ of the lattice 
$\Lambda(\alpha_1,\dots,\alpha_d)$ to the corresponding vector $\lambda'=(\lambda_1,
\dots,\lambda_d)$ with the same coordinates in $\Lambda(2,2,\dots,2)$
is a linear map which is one-to-one and length-shrinking.

This proves the Proposition if $\min(\alpha_1,\dots,\alpha_d)=2$.
In the remaining case, apply the Proposition to the lattice
$\Lambda'=\Lambda(\alpha_1-\kappa,\alpha_2-\kappa,\dots,\alpha_d-\kappa)$ where
$\kappa=\min(\alpha_1,\dots,\alpha_d)-2$ and observe that the obvious linear
map from $\Lambda(\alpha_1,\dots,\alpha_d)$ onto $\Lambda'$ is one-to-one
and diminuishes norms of non-zero elements at least by $\kappa$.
\hfill$\Box$

\begin{prop} Blocks of adjacents $2's$ in the sequence
$\alpha_1,\dots,\alpha_d$ are in bijection with irreducible root-sublattices
in $\Lambda(\alpha_1,\dots,\alpha_d)$.
\end{prop}

The easy proof is left to the reader.

\begin{thm} \label{thmlattice} Two lattices $\Lambda(\alpha_1,\dots,\alpha_d)$ and 
$\Lambda(b_1,\dots,b_d)$ of type T are isomorphic if and only if 
either $\alpha_i=b_i$ for $i=1,\dots,d$ or $\alpha_i=b_{d+1-i}$ for 
$i=1,\dots,d$.
\end{thm}

\begin{cor} The number of isomorphism classes
of Euclidean lattices (of arbitrary dimension) of type T,
determinant $N\geq 2$ and containing no elements of norm $1$
is given by $(\Phi(N)+\nu_N)/2$ where $\nu_N$ equals the number 
of elements of order $\leq 2$ in the multiplicative group 
$(\mathbb Z/N\mathbb Z)^*$.
\end{cor}

{\bf Idea for the proof of Theorem \ref{thmlattice}} We show by
induction on the dimension $d$ that metric properties of a lattice 
$\Lambda(\alpha_1,\dots,\alpha_d)$ of type T determine the integer sequence  
$\alpha_1,\dots,\alpha_d$ uniquely, up to reversion of the order.

If $\min(\alpha_1,\dots,\alpha_d)=2$, the root system $\mathcal R$ 
of $\Lambda(\alpha_1,\dots,\alpha_d)$ is determined by the set of blocks
of adjacent $2's$ in the sequence $\alpha_1,\alpha_2,\dots$. 
The sublattice $\tilde \Lambda$ orthogonal to the root
lattice generated by $\mathcal R$ corresponds to all 
terms $\alpha_i$ such that $\alpha_{i-1},\alpha_i,\alpha_{i+1}\geq 2$ and is an
orthogonal sum of lattices of type $T$. By induction 
on the rank, the lattice $\tilde \Lambda$ determines 
(up to reversion of the order)
disjoint subsequences $\alpha_{i_j},\alpha_{i_j+1},\alpha_{i_j+d_j}$
containing no $2$'s and not adjacent to $2$'s in 
$\alpha_1,\dots,\alpha_d$. Terms adjacent to irreducible root lattices
can essentially be recovered by considering minimal norms of sets of
vectors not contained and not orthogonal to such an irreducible 
root lattice. Orthogonality consideration allow then to glue 
all pieces together in an essentially unique way.

If $\min(\alpha_1,\dots,\alpha_d)>2$ one can replace the root lattice
by the lattice generated 
by the set of minimal vectors (corresponding to basis vectors
$\pm e_i$ such that $\alpha_i=\min(\alpha_1,\dots,\alpha_d)$) and proceed as above.
\hfill$\Box$

\begin{rem} Lattices of type T have poor
densities and are thus not interesting from the point of view of 
sphere-packings.
\end{rem}

\section{Polytopes}\label{sectpolyt}

Let $A$ be an integral symmetric matrix of size $d\times d$.
For every natural integer $N$, we denote by 
$\mathcal L_A(N)$ the set of integral diagonal matrices
of size $d\times d$
such that $A+D$ is positive definite of determinant $N$ for every
matrix $D\in\mathcal L_A(N)$. We denote by $\mathrm{Conv}({\mathcal L}_A(N))
\subset \mathbb R^d$ the convex hull of $\mathcal L_A(N)$.

\begin{thm} \label{thmpolyt}
For all $N$, the set $\mathrm{Conv}({\mathcal L}_A(N))$
is a polytope whose set of integral elements is contained in
$$\cup_{k=N}^\infty{\mathcal L}_A(k)$$
and whose vertices are given by the the integral elements $\mathcal L_A(N)$.
\end{thm}

{\bf Proof} The Brunn-Minkowski Theorem, see eg.
Theorem 6.2 in \cite{SY} states that
$$\det(A+\sum_{D\in\mathcal L_A(N)}\lambda_D D)^{1/d}\geq
\sum_{D\in\mathcal L_A(N)}\lambda_D \det(A+D)^{1/d}$$
if $\sum_{D\in\mathcal L_A(N)}\lambda_D=1$ with $\lambda_D$ positive real
numbers. This shows that we have $\det(A+D')\geq N$
for all $D'\in\mathrm{Conv}(\mathcal L_A(N))$ and the inequality
is strict if $D'\not\in \mathcal L_A(N)$.
Finiteness of the set $\mathcal L_A(N)$ follows from the fact 
that given two elements $(\alpha_1,\dots,\alpha_d)$ and
$(\beta_1,\dots,\beta_d)$ of $\mathcal L_A(N)$, there exists
indices $i,j\in\{1,\dots,d\}$ such that
$\alpha_i<\beta_i$ and $\alpha_j>\beta_j$.\hfill$\Box$

The construction of the polytope $\mathrm{Conv}(\mathcal L_A(N))$
has the following properties, given without proofs. (They are
straightforward.)

The diagonal part of $A$ is essentially irrelevant:
If $\tilde A=A+D$ where $D$ is integral diagonal and $A$ is
integral symmetric, then $\mathcal L_{\tilde A}(N)=
\mathcal L_A(N)-D$. We suppose henceforth that $A$ has
zeroes along the diagonal.

Automorphisms of $A$ (conjugations by signed permutation matrices 
commuting with $A$)
act on the polytopes $\mathrm{Conv}(
\mathcal L_A(N))$ in the obvious way by permuting the coefficients
according to the underlying ordinary permutation matrix.

If $A$ is a ``direct sum'' $A=A_1+A_2$ of disjoint diagonal blocks, then 
$$\mathrm{Conv}(\mathcal L_A(1))=\mathrm{Conv}(\mathcal L_{A_1}(1))
\times \mathrm{Conv}(\mathcal L_{A_1}(1))\ .$$
More generally,
$$\mathcal L_A(N)=\cup_{l\vert N}\mathcal L_{A_1}(l)\times
\mathcal L_{A_2}(N/l)$$
where the union is over all natural integers divising $N$.

The following construction yields a map from 
the permutation group on the $d$ indices into the set $\mathcal L_A(1)$:
Given a permutation matrix $\sigma$ of size $d\times d$, we
consider the unique integral diagonal matrix $D=D_\pi$
such that the submatrix formed by the first $k$ rows and columns of 
$\sigma^{-1}(A+D_\sigma)\sigma$ is unimodular for every integer
$k$ in $\{1,\dots,d\}$. These special vertices are invariant under
automorphisms of $A$ (and their convex hull is a polytope 
contained in $\mathrm{Conv}(\mathcal L_A(1))$ which is
invariant under automorphisms of $A$).

If $A$ is an adjacency matrix of a simple finite undirected graph
which is connected,
the above construction can be restricted to permutations
such that the first $k$ vertices form connected subgraphs for 
all $k$. The subset of these vertices
is still invariant under automorphisms. There exists slight
generalisations of this construction. 

{\bf Examples} The easiest case is given if the matrix $A$
is the zero-matrix of size $d\times d$. 
For $d\geq 2$, the polytope $\mathrm{Conv}(\mathcal L_A(N))$
is a $(d-1)-$dimensional simplex if and only if $N$ is a 
prime number. More interesting cases occur if $N$ is highly composite. 

A probably interesting example is given
by considering the symmetric matrix $A$ of size $d\times d$
with coefficients $A_{i,j}=1$ if $\vert i-j\vert=0$ and 
$A_{i,j}=0$ otherwise. By \cite{LN}, the set $\mathcal L_A(1)$ contains 
then exactly $C_d={2d\choose d}/(d+1)$ elements forming the vertex set 
of an integral polytope $P_d$ of dimension $d$ if $d\geq 3$.
The polytope $P_3$ for instance is the convex hull of the five
vertices
$$\begin{array}{l}
v_1=(1,2,2)\\
v_2=(1,3,1)\\
v_3=(2,1,3)\\
v_4=(2,2,1)\\
v_5=(3,1,2)\end{array}$$
Its faces are defined by equalities of the five inequalities
$$\begin{array}{ll}
y+z\leq 4&\{v_1,v_2,v_3\}\\
x+y+z\geq 5&\{v_1,v_2,v_4\}\\
x+2y+z\geq 7&\{v_1,v_3,v_4,v_5\}\\
2x+3y+2z\leq 13&\{v_2,v_3,v_5\}\\
x+y\leq 4&\{v_2,v_4,v_5\}\end{array}$$
(with vertex-sets indicated for each face)
and the polytope $P_3$ is thus a pyramid with summit $v_2$ over
the planar polygone formed by the four vertices $v_1,v_3,v_4,v_5$.
The face $\{v_1,v_2,v_4\}$ contains all vertices associated
to permutation matrices of size $3\times 3$. (This can easily be
generalised to an arbitrary dimension $d$: 
Such special vertices are all of the form
$(1,2,2,2,\dots,2,1)+e_i$ where $e_i=(0,\dots,0,1,0,\dots,0)$
denotes the $i-$th basis vector and form a unique simplicial
face of the of polytope $\mathrm{Conv}(\mathcal L_A(1))$.)

It would perhaps be interesting to understand the polytopes $P_d$
in general,
in particular eventual connections with so-called
Associahedrons (or Stasheff polytopes) whose vertices are also enumerated
by Catalan numbers and whose dimensions are in general one less.

\noindent Roland BACHER

\noindent INSTITUT FOURIER

\noindent Laboratoire de Math\'ematiques

\noindent UMR 5582 (UJF-CNRS)

\noindent BP 74

\noindent 38402 St Martin d'H\`eres Cedex (France)
\medskip

\noindent e-mail: Roland.Bacher@ujf-grenoble.fr

\end{document}